\documentclass[10pt, oneside, draft]{amsart}

        \usepackage {amssymb,amsthm,amsxtra,amscd,amsmath,amsfonts}
        
        \newcommand{\K}{\ensuremath{\Bbbk}}
        \newcommand{\Q}{\ensuremath{\mathbb{Q}}}
        \newcommand{\Si}{\ensuremath{\widetilde{\sigma}}}
        \newcommand{\Z}{\ensuremath{\mathbb{Z}}}
        
        \newcommand{\CA}{\ensuremath{\mathcal{A}}}
        
        \newcommand{\CZ}{\ensuremath{\textup{Z}}}
        
        \newcommand{\degr}[1]{\textup{deg}(#1)}

        \newcommand{\grhom}{\textup{grHom}_\K(R, R)}
        \newcommand{\op}{\textup{o}}
        \newcommand{\cross}{\times}
        \newcommand{\tensor}{\bigotimes}
        \newcommand{\mb}[1]{\mathbf{#1}}

        \theoremstyle{plain}
                \newtheorem{theorem}{Theorem}[subsection]
                \newtheorem{lemma}{Lemma}[subsection]
                \newtheorem{corollary}{Corollary}[subsection]
                \newtheorem{proposition}{Proposition}[subsection]
                \newtheorem{remark}{Remark}[subsection]
                \newtheorem*{definition}{Definition}
                \newtheorem*{theorem*}{Theorem}
        \theoremstyle{definition}
                \newtheorem*{note}{Note}
        \numberwithin{equation}{subsection}
        
        \newcommand{\ignore}[1]{}
        \newcommand{\term}[1]{\textit{#1}}
        \newcommand{\mynote}[1]{}

\begin{document}
\setcounter{section}{-1}


        \title[$q$-differential operators]{Quantum differential operators on 
                        $\K[x]$}
        \date{\today}
        \author[Iyer]{Uma N. Iyer} 
        \email{uiyer@mri.ernet.in}
        \author[McCune]{Timothy C. McCune} 
        \email{tim@mri.ernet.in}
        \address{Harish-Chandra Research Institute, Chhatnag Road, Jhusi,\\
                  Allahabad, U.P. 211 019\\
                  India}

\begin{abstract}
Following the definition given in \cite{LR1}, we compute the ring of
quantum differential operators on the polynomial ring in 1 variable.  
We further study this ring.
\end{abstract}
\maketitle              

\section{Introduction}\label{sect:introduction}
Quantum differential operators were defined by V.A.Lunts and A.L.Rosenberg
in their work \cite{LR1}, as part of their project on Localization for Quantum
groups (\cite{LR2}).
Quantum differential operators are defined on graded rings, graded by an 
abelian group.  
Let $\K$ be a field and $R$ an associative \K -algebra, graded by an abelian 
group $\Gamma$.  
The ring of quantum differential operators (or $q$-differential operators),
denoted by $D_q(R)$ is defined for a fixed
bicharacter $\beta :\Gamma \times \Gamma \to \K ^*$.
For each $a\in \Gamma$, there is an automorphism $\sigma _a$ of $R$
given by $\sigma _a (r_b) = \beta (a,b) r_b$, where degree of $r_b = b$.
Turns out that each $\sigma _a$ and  the corresponding 
left-$\sigma _a$-derivation is a quantum differential operator.
(This is shown in section \ref{sect:preliminaries}).

In this paper we compute the ring of
$q$-differential operators on a polynomial ring
in 1 variable.  We also show a relationship between this ring and the
quantum group on $sl_2$.
  
In the first section
we cover the preliminaries required in the rest of the paper

In the second section,
we consider the polynomial ring on one variable, graded by $\mathbb{Z}$.
Let  $q$ be transcendental over $\Q$ and $\K$ be a field containing $\Q(q)$. 
We let $R = \K [x]$, graded by the group of integers as $\mathit{deg} (x) =1$.
The bicharacter $\beta$ is defined as $\beta (n,m) = q^{nm}$.
Using this set up, we show in Theorem \ref{T:qdops} 
that $D_q(R)$ is generated over \K by
$x,\partial ,\partial ^{\beta}, \partial ^{\beta ^{-1}}$, where
$\partial$ is the usual derivation, ($\partial (x^n) = nx^{n-1}$),
$\partial ^{\beta}$ is the left $\sigma _1$-derivation
($\partial ^{\beta}(x^n) = (1+q+\cdots + q^{n-1})x^{n-1}$) and
$\partial ^{\beta ^{-1}}$ is the left $\sigma _{-1}$-derivation
($\partial ^{\beta}(x^n) = (1+\frac{1}{q}+\cdots + \frac{1}{q^{n-1}})
x^{n-1}$).

Here, we also show that $D_q(R)$ is a simple ring (\ref{T:Dqsimple}), a result
similar to that of usual differential operators on $R$ (characteristic 0).

In section \ref{sect:intrinsic}, our main goal is to write $D_q(R)$ in terms
of generators and relations, given in theorem \ref{T:intrinsic};
which is,
\begin{theorem*}
        The ring $D_q$ of $q$-differential operators on $R$ is 
        the $\K$-algebra generated by $x$, $\partial^{\beta^{-1}}$, 
        $\partial^{\beta^0} (= \partial )$, and  
        $\partial^{\beta^1}$ subject to the relations
        \begin{align*}
                \partial^{\beta^{a}} x - q^{a} x \partial^{\beta^{a}}
                        &= 1, &
                \partial^{\beta^{a}} x \partial^{\beta^{b}} 
                        &= \partial^{\beta^{b}} x  \partial^{\beta^{a}}, &
                        &\text{and} &
                \partial^{\beta^{-1}} - q  \partial^{\beta} 
                        &= (1 - q)\partial^{\beta^{-1}} x  \partial^{\beta}.
        \end{align*}
\end{theorem*}
Other results of independent interest are also proved in this section.
If $\tau = x \partial$ and $\sigma = \sigma _1$, and $(D_q)_0$ denotes the
$q$-differential operators of degree 0, then
lemma \ref{L:D0domain} shows that $(D_q)_0$ is localization of the polynomial
 ring
$\K [\tau ,\sigma]$ at the element $\sigma$; that is
$(D_q)_0 = \K [\tau, \sigma ,\sigma ^{-1}]$.
A corollary (\ref{C:Dqdomain}) to this lemma is 
the fact that $D_q(R)$ is a domain (again similar to the result in 
usual differential operators).

In section \ref{sect:nvariables}, we generalize the computation of $D_q(R)$
to the case when $R$ is a polynomial ring in $n$  variables, graded
by $\mathbb{Z}^n$, and the 
bicharacter is product of $n$-bicharacters of the above type.
That is, we let $\K$ let be a field containing $\mathbb{Q}$ and $n$ 
transcendental elements $q_1,q_2,\cdots ,q_n$, and
define 
\[
\beta ((a_1,a_2,\cdots ,a_n),(b_1,b_2,\cdots ,b_n)) =
          q_1^{a_1b_1}q_2^{a_2b_2}\cdots q_n^{a_nb_n}.
\]
Here, again as in the case of usual differential operators, one can define
for each $i,1\leq i \leq n$, the operators $\partial ^{\beta ^{-1}}_i,
\partial _i,\partial ^{\beta ^1}_i$. We show  in theorem \ref{T:nvariables}
that
these operators along with $R$ generate $D_q(R)$.

In the last section we explain relationship between quantum group on 
$sl_2(\K)$,
denoted by $U_q(sl_2)$,
and the global $q$-differential operators on $\mathbb{P}^1$.
Here, we fix $\Gamma = \mathbb{Z}$ and $\beta (n,m) = q^{nm}$.
The $q$-differential operators $D_q(\K [x])$ and $D_q(\K [x^{-1}])$
extend to $q$-differential operators on $\K [x,x^{-1}]$ by theorem 3.2.2. of
\cite{LR1}.
So, we let $\Gamma _q(\mathbb{P}^1)$ to be the kernel of
\begin{align*}
D_q(\K [x]) \bigoplus D_q(\K [x^{-1}]) &\to D_q(\K [x,x^{-1}])\\
(\varphi _1,\varphi _2) &\mapsto \varphi _1 - \varphi _2.
\end{align*}
There is a homomorphism $\eta : U_q(sl_2) \to \Gamma _q (\mathbb{P}^1)$,
which generalizes the homomorphism of the enveloping algebra of $sl_2$ to
the ring of global sections of usual 
differential operators on $\mathbb{P}^1$.
But unlike the enveloping algebra case (characteristic 0), the
homomorphism $\eta$ is not a surjection. So, next we replace \K
by the ring $\mathcal{A} = \mathbb{Q}[q,q^{-1}]_{(q-1)}$ a local ring where
$q-1$ is not invertible.  We consider the corresponding rings $U_q(sl_2)$
and $\Gamma _q$ with $\mathcal{A}$ as the base ring
and consider their respective inverse limits with respect to the ideal
$(q-1)$.  The homomorphism $\eta$ gives rise to a homomorphism of 
inverse limits, denoted by $\hat{\eta}$.  We show that
$\hat{\eta}$ is a surjection in theorem \ref{T:inverselimits}.
We would like to point out that some of the formulae derived here have
already been derived in the literature of Mathematical Physics, see for 
example \cite{GP}.

We thank Professor Valery Lunts for suggesting this problem especially
for  his ideas
regarding the relationship with quantum group on $sl_2$.


\section{Preliminaries}\label{sect:preliminaries}

Throughout this paper,
let $q$ be transcendental over $\Q$ and $\K$ be a field containing $\Q(q)$. 
Let $R$ be a $\K$-algebra which is graded by an abelian group $\Gamma$.
That is, $R = \oplus_{\mb{a} \in \Gamma} R_{\mb{a}}$ and
$R_\mb{a} R_\mb{b} \subseteq R_{\mb{a}+\mb{b}}$.

A $\K$-linear endomorphism $\varphi$ of $R$ is called
\term{homogeneous of degree $\mb{a} \in \Gamma$}
if and only if for every $\mb{b} \in \Gamma$, 
$\varphi(R_\mb{b}) \subseteq R_{\mb{a}+\mb{b}}$.
We will say $\varphi$ is a \term{graded endomorphism} if 
and only if it is the sum of homogeneous endomorphisms.  
The $\K$-vector space of all graded endomorphisms of $R$ 
will be denoted $\grhom$. 

For each $r \in R$, we have two endomorphisms of $R$, 
left multiplication by $r$ and right multiplication by $r$.
We will denote these by $\lambda_r$ and $\rho_r$ respectively.  
When $r$ is homogeneous of degree $\mb{a}$, so are $\lambda_r$ and 
$\rho_r$.  Hence, for any $r$, the maps $\lambda_r$ and 
$\rho_r$ are graded endomorphisms.
The homomorphism $r \mapsto \lambda_r$ 
embeds $R$ into $\grhom$.  This enables us to define right and 
left actions of $R$ on $\grhom$ by 
$r \cdot \varphi \cdot s = \lambda_r \varphi \lambda_s$.
Since a right $R$ action is equivalent to a left action of the 
oposite ring $R^{\op}$, we have a left action of the 
enveloping algebra $R^{\textup{e}} = R \tensor_\K R^{\op}$
on $\grhom$. 
\mynote{Similarly, the homomorphism, $r \mapsto \rho_r$ embeds 
        $R^{\op}$, into $\grhom$ yielding another canonical 
        action of $R^{\textup{e}}$ on $\grhom$, but we shall only 
        make use of the former action in this paper.}

A bicharacter on $\Gamma$ is a function $\beta$ from 
$\Gamma \cross \Gamma$ to $\K^{\cross}$ such that, in each 
argument, $\beta$ is a group homomorphism.  To each bicharacter 
we can associate a family of automorphisms of $\grhom$ as follows.  
For each $\mb{a} \in \Gamma$, define $\Si_\mb{a}$
so that for any homogeneous $\varphi$ of degree $\mb{b}$, 
$\Si_{\mb{a}}(\varphi) = \beta(\mb{a}, \mb{b}) \varphi$.
This extends linearly to all of $\grhom$.

Now with each $\Si_\mb{a}$ we can define a pairing in $\grhom$.
Let $\varphi$ and $\psi$ be graded endomorphisms. Define the 
\term{$\mb{a}$-twisted bracket} to be 
\[      
        [\varphi, \psi]_\mb{a} = 
                \varphi \psi - \Si_\mb{a}(\psi) \varphi.
\]
The reader is cationed that this will be a Lie bracket 
only if $\mb{a}$ is the identity $\mb{0}$. In this case, the 
$\mb{0}$-twisted bracket will simply be denoted by 
$[ \cdot, \cdot]$.

In a similar way, we can define a family of endomorphisms of $R$.
If $r \in R$ is homogeneous of degree $\mb{b}$, then define
$\sigma_{\mb{a}}(r) = \beta(\mb{a},\mb{b})r$.  
This can be extended linearly to all of $R$ making $\sigma_\mb{a}$
an homogeneous endomorphism of degree $0$.  Then for any $r \in R$, 
$\Si_{\mb{a}}(\lambda_r) = \lambda_{\sigma_{\mb{a}}(r)}$.
More generally, for any graded endomorphism $\varphi$ we have
\begin{equation}\label{prelim:1} 
        \Si_{\mb{a}}(\varphi)\sigma_{\mb{a}} = \sigma_\mb{a} \varphi.
\end{equation}
\mynote{This is telling of the fundamental r\^{o}le
        the endomorphisms $\sigma_{\mb{a}}$ play in the structure 
        of the quantum differential operators on $R$.}

\subsection{Quantum differential operators}
Following \cite{LR1}, we define $D_q^0(R)$, the \term{$q$-centre} of
$\grhom$, to be the smallest $R^{\textup{e}}$ submodule of $\grhom$
containing the $\K$-space
\mynote{
        Following \cite{LR1}, we define the \term{$q$-centre} of 
        $\grhom$ to be the $\K$-space}
\begin{align*}
        \CZ_q^0 = \K \langle \text{ homogeneous $\varphi$} \mid 
                & \text{ there is some $\mb{a} \in \Gamma$ 
                        such that for}\\ 
                & \text{ any $r \in R$, 
                        $[\varphi, \lambda_{r}]_\mb{a} =0$} \rangle .
\end{align*}
\mynote{Then $D_q^0(R)$ is defined to be the smallest 
        $R^{\textup{e}}$-submodule of $\grhom$ containing $\CZ_q^0$.}
Suppose that we have defined already the $R^{\textup{e}}$-submodule
$D_q^n(R)$ of $\grhom$.  Then the $R^{\textup{e}}$-module $D_q^{n+1}(R)$ is 
defined to be the smallest
$R^{\textup{e}}$-submodule containing the $\K$-space
\begin{align*}
        \CZ_q^{n+1} = \K \langle \text{ homogeneous $\varphi$} \mid 
                & \text{ there is some $\mb{a} \in \Gamma$ such 
                        that for any}\\ 
                & \text{ $r \in R$, 
                        $[\varphi, \lambda_{r}]_\mb{a} \in D_q^n(R)$}\rangle .
\end{align*}
We say $\varphi$ is a \term{$q$-differential operator of order $n$} if 
and only if $\varphi \in D_q^n(R)$.  
Finally, the set of all $q$-differential operators is defined to be 
$D_q(R) = \bigcup D_q^n(R)$.

\subsection{The $q$-centre of $D_q(R)$}
As is the case with differential operators over a commutative ring, 
the composite of a $q$-differential operator of order $n$ with one 
of order $m$ is a $q$-differential operator of order $m+n$ 
\cite[3.1.8]{LR1}.
This shows that $D_q^0(R)$ is a ring and each $D_q^n(R)$ is a 
$D_q^0(R)$-module.  

The structure of $D_q^0(R)$ is apparent:
\begin{lemma}\label{L:Dq0}
        The ring $D^0_q(R)$ is a $\Gamma$-graded $\K$-algebra 
        generated by 
        \[
                \{\lambda_r \rho_s \sigma_\mb{a} \mid 
                \text{ $\mb{a} \in \Gamma$ and $r,s \in R$ }\} .
        \]
\end{lemma}
\begin{proof}
        First, we will show that $\CZ_q^0$ is generated as a $\K$-module
        by 
        \[
                \{ \rho_r \sigma_{\mb{a}}\mid \text{ $r$ is homogeneous } \}.
        \] 

        Let $\varphi$ be a homogeneous endomorphism in the $q$-centre of
        $\grhom$.  
        Then,  $[\varphi , \lambda_s]_\mb{a} = 0$ 
        for some $\mb{a} \in \Gamma$ and all $s$ in $R$.  
        Since $\varphi$ is homogeneous, $\varphi(1)$ is some 
        homogeneous $r$ of the same degree as $\varphi$.  
        Then $\varphi (s) = \varphi \lambda_s(1)
        = \lambda_{\sigma_\mb{a}(s)}\varphi(1)
        = \sigma_\mb{a}(s) r$.
        Since $\varphi$ has the same value on $s$ as 
        $\rho_{r}\sigma_\mb{a}$, we have 
        $\varphi = \rho_{r} \sigma_\mb{a}$.

        To prove $\CZ_q^0$ contains all of the prescribed operators, 
        let $r$ be any homogeneous element 
        and $\mb{a}$ any element of $\Gamma$.
        Since left multiplication and right multiplication commute, 
        for any $s$ we have 
        $ [\rho_{r} \sigma_\mb{a}, \lambda_s]_\mb{a} = 
        \rho_r[\sigma_\mb{a}, \lambda_s]_\mb{a}$.
        Furthermore, by (\ref{prelim:1}) we have 
        $[\sigma_\mb{a}, \lambda_s]_\mb{a} = 0$.
        Hence, $\rho_{r} \sigma_\mb{a}$ is in the $q$-centre of $\grhom$.

        Now, by its definition $D_q^0(R)$ is the $\K$-subspace of $\grhom$ 
        spanned by 
        $\{ \lambda_r \rho_s \sigma_\mb{a} \lambda_t \}$.  However, 
        $\lambda_t$ commutes with $\sigma_\mb{a}$ up to multiplication 
        by a scalar, and commutes with $\rho_s$ without incident, so 
        $D_q^0(R)$ can be generated by the given operators as required.
\end{proof}

\begin{corollary}\label{C:1}
        The module $D_q^{n+1}(R)$ of $q$-differential operators of order 
        $n+1$ are generated over $D_q^0(R)$ by 
        \[
                {\CZ'_q}^{n+1} = \{ \text{ homogeneous $\varphi$} \mid 
                        \text{ such that for any $r \in R, 
                         [\varphi, \lambda_{r}] \in D_q^n(R)$} \}.
        \]
\end{corollary}
\begin{proof}  
        It is clear that ${\CZ'_q}^{n+1}$ is contained in $\CZ_q^{n+1}$ so
        we only need to show that  
        every  $\varphi \in \CZ_q^{n+1}$ is in the $D_q^0(R)$ span of 
        ${\CZ'_q}^{n+1}$.  To this end, suppose for some $a \in \Gamma$
        and all $r \in R$ that $[\varphi, \lambda_r]_\mb{a} \in D_q^n(R)$.
        Then we have for any $b \in \Gamma$,
        \begin{align*}
                [\varphi \sigma_\mb{b}, \lambda_r]_{\mb{a}+\mb{b}} = &
                \varphi \sigma_\mb{b} \lambda_r 
                     - \Si_{\mb{a}+\mb{b}} (\lambda_r) \varphi\sigma_\mb{b}\\
                =  & \varphi \Si_\mb{b}(\lambda_r)\sigma_\mb{b}
                     - \Si_\mb{a}(\Si_\mb{b}(\lambda_r))\varphi \sigma_\mb{b}\\
                = & [\varphi, \Si_\mb{b}(\lambda_r)]_\mb{a} \sigma_\mb{b}.
        \end{align*}
        Since $\sigma_\mb{b}$ is an automorphism of $R$, 
        $[\varphi, \Si_\mb{b}(\lambda_r)]_\mb{a} \in D_q^n(R)$.
        Hence $[\varphi \sigma_\mb{b}, \lambda_r]_{\mb{a}+\mb{b}} 
        \in D_q^n(R)$.  Since $\Gamma$ is a group, we may put 
        $\mb{b} = -\mb{a}$.  Then $[\varphi \sigma_{-\mb{a}}, \lambda_r]
        \in D_q^n(R)$ for any $r \in R$.  Hence $\varphi \sigma_{-\mb{a}}
        \in {\CZ'_q}^{n+1}$.
\end{proof}

\begin{remark}
        Typically, a bicharacter $\beta$ is used to define an 
        associative algebra structure on the tensor product
        of two graded algebras.
        Let $A = \oplus_{\mb{a}\in\Gamma} A_\mb{a}$ be a 
        $\Gamma$-graded $\K$-algebra, and  define $A_{\Gamma}$ to be the 
        algebra whose underlying set of elements is $A \tensor \K[\Gamma]$,
        and whose multiplicative structure is given as follows:
        If we define the grading automorphisms $\sigma_\mb{a}$ as above, 
        then define 
        \[
         r \otimes \mb{a} \cdot s \otimes \mb{b} 
                = r (\sigma_{\mb{a}}(s))
                         \otimes (\mb{a}+\mb{b}).
        \]
        Requiring $\beta$ to be a group homomorphism in each variable
        ensures that $A_\Gamma$ is associative.  We call this the 
        \term{crossed-product algebra} determined by $\beta$.
        
        It is not surprising in light of this that $D_q^0(R)$ has 
        the structure of a crossed-product algebra.  Indeed, if 
        we put $A = R\tensor_{Z(R)} R^{\op}$ where $Z(R)$ is 
        the center of $R$, then $D_q^0(R) = A_\Gamma$.
\end{remark}

\subsection{Left $\beta$-derivations}
        The definition of a derivation of a ring is quite standard and 
        can be applied to commutative and noncommutative rings alike.  
        A derivation is always a differential operator, and, in the 
        best of circumstances, the derivations generate all other 
        differential operators.
        However, since the $q$-differential operators form a broader 
        class of graded endomorphisms, even in the best of circumstances, 
        we will require more than just the derivations to generate
        all of the $q$-differential operators.

\begin{definition}
        A left $\beta$-derivation is a graded endomorphism $\varphi$ 
        of $R$ such that $\varphi$ obeys a ``twisted'' Leibniz's Rule.
        That is, for some $a \in \Gamma$ and any $r, s \in R$, 
        \[      \varphi(rs) = \varphi(r)s + \sigma_\mb{a}(r)\varphi(s).
        \]
\end{definition}
        Equivalently, we could define a left $\beta$-derivation 
        to be an endomorphism $\varphi$ such that for some
        $\mb{a} \in \Gamma$ and any $r \in R$ we have 
        $[\varphi, \lambda_r]_\mb{a} = \lambda_{\varphi(r)}$.
        From this description, it becomes clear why let $\beta$-derivations
        are $q$-differential operators, and why the analogously defined
        \textit{right} $\beta$-derivations, 
        those $\varphi$ for which there is an $\mb{a} \in \Gamma$
        such that for all $r \in R$ we have
        $[\varphi, \rho_r]_\mb{a}= \rho_{\varphi(r)}$,
        bear little importance to us.


\section{The $q$-Differential Operators on a Polynomial Ring in 1 Variable}
        \label{S:1var}

        Let $R = \K[x]$, and, to simplify notation, let $D^n_q = D^n_q(R)$ 
        and $D_q = D_q(R)$.
        The ring $R$ has a $\Z$-grading given by $\degr{x^a} = a$.
        Denote by $R_a$ the $a$-th graded part $\K \cdot x^a$.
        Let $\beta :\Z \times \Z \to \K^{\cross}$ be a bicharacter.
        Since $\beta (n,m) = (\beta(1,1))^{nm}$,
        we let $\beta (n,m) = q^{nm}$.

        Since $R$ is commutative, $R^{\op} = R$, and 
        $D_q^0 = \K\langle r \sigma_a \mid 
                \text{ $r \in R$ and $a \in \Z$ } \rangle$.
        In this section, we shall give an explicit description
        of $D_q$.

        We begin with two lemmas which will help us recognize 
        $q$-differential operators of order $n$.
\begin{lemma}\label{L:gens}
        If $[\varphi, x] \in D_q^{n}$ then $\varphi \in \CZ_q^{'n+1}$.
\end{lemma}
\begin{proof}
        Since $[ \cdot, \cdot]$ is a bracket, $[\varphi, \cdot]$ is 
        a derivation.  Hence, for any $m >0$, we have
        \[ 
                [\varphi, x^m] = [\varphi, x] x^{m-1} + x [\varphi, x^{m-1}].
        \]
        Inductively, we have that if $[\varphi, x] \in D_q^{n}$ 
        then $[\varphi, x^{m-1}] \in D_q^{n}$.  Hence $[\varphi, x^m]$
        is in the $R$-span of $D_q^{n}$, and so it is in $D_q^{n}$.
        It follows that for any $r \in R$, $[\varphi, r] \in D_q^{n}$.
        Hence $\varphi \in \CZ_q^{'n+1}$.
\end{proof}
        This means we do not have to test $[\varphi, \cdot]$ against 
        all elements of $R$ in order to determine whether or not
        $\varphi$ is in $\CZ_q^{'n+1}$.  For general rings $R$, we have 
        a similar result which we will not prove: 
        $\varphi \in \CZ_q^{'n+1}$ if and only if 
        $[\varphi, \lambda_r] \in D_q^{n}$ 
        for all generators $r$ of $R$ over $\K$.   

\begin{lemma}\label{L:comm}
        Given $\varphi \in D_q(R)$, there exists integers
        $n,a_1,a_2,\cdots , a_n$ such that $n > 0$ and  
        $[\cdots[ [\varphi , x]_{a_1}, x]_{a_2}, 
         \cdots ,x]_{a_n} = 0$.
\end{lemma}
\begin{proof}
Indeed, if $\varphi \in D^k_q(R)$, then 
\[
\varphi = \varphi _1 +\varphi _2 +\cdots +\varphi _l
\]
 for homogeneous
$\varphi _i \in \CZ_q^{k}$.  Corresponding to each $\varphi _i$, there
exists $b_i \in \Z$ such that $[\varphi_i, x]_{b_i} 
\in D^{k-1}_q(R)$.  
Hence, 
\[
[\cdots [[\varphi , x]_{b_1}, x]_{b_2},
\cdots , x]_{b_l} \in D^{k-1}_q(R).
\]
Now induction completes the lemma.
\end{proof}

        Let us now describe the left $\beta$-derivations on $R$.
\begin{definition}
        For each $a \in \Z$, define $\partial^{\beta^a} \in 
        \grhom$ as
        \[
                \partial^{\beta^a}(x^b) = (1 + q^a +q^{2a} +
                \cdots + q^{a(b-1)})x^{b-1}.
        \]
        We denote $\partial^{\beta^1}$ by just $\partial^{\beta}$.
\end{definition}
        Simply by comparing values on $x^m$, one is led to 
        the following:
\begin{note}\label{N:nto1}\hfill
\begin{enumerate}
        \item   For any positive integer $a$, 
                \[ 
                \partial ^{\beta ^a} = \left( \frac{1-q}{1-q^a} \right)
                \partial ^{\beta}[1+\sigma_1 +\cdots +\sigma_{a-1}].
                \]
        \item   When $a = 0$, $\partial^{\beta^a}$ is $\partial$, the usual 
                derivative defined by $\partial(x^b) = bx^{b-1}$.
        \item   For any positive integer $a$,
                \[
                \partial ^{\beta ^{-a}} = \left( \frac{1-q}{1-q^{-a}} \right)
                \partial ^{\beta^{-1}}[1+\sigma_{-1} +\cdots +\sigma_{1-a}].
                \]
        \item   For any integer $a$, 
                \begin{equation}\label{tim:3}
                        \partial^{\beta ^a} 
                                = \sigma_{a}\partial^{\beta^{-a}}.
                \end{equation}
\end{enumerate}
\end{note}

\begin{lemma}\label{L:dinD1}
        For any $a \in \Z$, the operator $\partial^{\beta^a}$ is
        in $D_q^1$.
\end{lemma}
\begin{proof}
        Let $\eta = [\partial^{\beta^a}, x]$.  By  
        Lemma \ref{L:gens}, it is enough to show $\eta \in D_q^0$.
        For any $m \geq 0$, 
        \begin{align*}
                \eta(x^m) 
                        = & \partial^{\beta^a}(x^{m+1}) 
                           - x \partial^{\beta^a}(x^m)\\
                        = & (1 + q^a + \cdots + q^{am})x^m
                           - x (1 + q^a + \cdots + q^{a(m-1)})x^{m-1}\\
                        = & q^{am} x^m.
        \end{align*}
        Thus $\eta(r) = \sigma_a(r)$ for any $r \in R$.
        Hence $\eta = \sigma_a \in D_q^0$, as required. 
\end{proof}

        We will make use of the following notation.
        For each positive integer $n$ and each multi-index
        $I = (a_1,a_2,\cdots ,a_n)$, set
        \begin{equation}\label{E:multiindex}
                \partial ^{\beta ^I} =
                \partial ^{\beta ^{a_1}}\partial ^{\beta ^{a_2}} \cdots
                \partial ^{\beta ^{a_n}}.
        \end{equation}
        Here we say $|I| = n$.

\begin{corollary}\label{C:2}
        If $P = \partial^{\beta^I}$ and $n \geq |I|$, then $P \in D^n_q$.
\end{corollary}
\begin{proof}
        This follows from the fact that $P \in (D_q^1)^n \subset D_q^n$.
\end{proof}

%
%

        We will prove the following theorem:
\begin{theorem}\label{T:qdops}
        The ring $D_q$ of $q$-differential operators on $R$ is 
        generated as a $\K$-algebra by $R$ and the set
        $\{ \partial^{\beta^{-1}}, \partial , \partial ^{\beta} \}$.
\end{theorem}
\mynote{To say that $D_q$ is an $R$-algebra is to say that $R$ is in the 
        center of $D_q$.  This is not true in general, so $D_q$ is just
        a $\K$-algebra and $R$-module}
        The proof is somewhat involved, so we have consigned the more 
        technical parts of this theorem to a Lemma.
        We will use the following identity which can be easily checked by 
        expanding the $a$-twisted bracket $[\cdot, \cdot]_a$ according to 
        its definition:
        \begin{equation}\label{tim:2}
                [[\varphi, \psi]_a, x] 
                        = [[\varphi, x], \psi]_a +
                            \varphi [\psi, x] - [\Si_a(\psi), x] \varphi.
        \end{equation}

\begin{lemma}\label{L:intmono}
        If $P$ is a $q$-differential operator comprised of monomials in 
        $\{\partial^{\beta^{a}} \mid a \in \Z \}$ and $b$ is any integer
        then there is a $q$-differential
        operator $Q$, also comprised of monomials in the
        $\partial^{\beta^{a}}$'s with coefficients in $\K$, 
        such that $[Q,x] = P \sigma_b$.  Thus, for $f\in R$, we have
        $[fQ,x] = fP\sigma _b$.
\end{lemma}
\begin{proof}
        It is enough to prove the lemma when $P$ is a monomial, for 
        if the lemma holds for monomials $m_i$, $i = 1, \ldots, k$, and 
        $P = m_1 + m _2 + \cdots + m_k$, then we can find differential 
        operators $M_i$ of the prescribed form
        such that $[M_i, x] = m_i \sigma_b$.  Putting
        $Q = M_1 + M_2 + \cdots + M_k$ yields the desired result.

        Now, suppose $P = d_{n}d_{n-1} \cdots d_{2} d_{1}$ where 
        each $d_i = \partial^{\beta^{a_i}}$.  
        We will show by induction on $n$ that there is a homogeneous
        $q$-differential operator $Q$ of degree $n+1$, comprised of 
        monomials
        in the $\partial^{\beta^{a}}$'s such that $[Q,x]= P \sigma_b$.

        When $n=0$ we have only one possibility:  $P = 1$.
        Put $Q = \partial^{\beta^{b}}$.
        Since $[Q, x] = \sigma_b = P \sigma_b$, this proves the base case.

        Before proceeding, we must address one technical point.
        We shall assume that if 
        $b = -\sum_{i=1}^n a_i$
        then 
        every $a_i = 0$.  To see that no generality is lost, 
        suppose that $b = -\sum a_i$ and  $a_n\neq 0$. 
        Put $P' = \partial^{\beta^{-a_n}} d_{n-1}\cdots d_2 d_1$.
        By (\ref{tim:3})
        $\partial^{\beta^{a_n}} = \sigma_{a_n} \partial^{\beta^{-a_n}}$.
        Thus we have 
        \[
                P \sigma_b 
                = \sigma_{a_n} \partial^{\beta^{-a_n}} d_{n-1} \cdots d_1 
                        \sigma_b
                =  c P' \sigma_{a_n + b}
        \] 
        for some nonzero constant $c$.  
        If $[Q', x] = P' \sigma_{a_n + b}$ 
        then $[cQ' , x] = P \sigma_b$. 
        However, since $a_n >0$, we have 
        $-(-a_n + a_{n-1} + \cdots +a_1) = 2 a_n -\sum a_i \neq  a_n +b$.

        Now suppose the statement holds for monomials of degree $n-1$.
        Let 
        \begin{align*}
                t_1 = & d_n d_{n-1} \cdots d_2,\\
                t_2 = & d_1 \; d_n \; \cdots \; d_3,\\
                \vdots & \\ 
                t_n = & d_{n-1} d_{n-2} \cdots d_1.\\
        \end{align*}
        That is $t_i d_i$ is the monomial obtained from $P$
        by a cyclic permutation of its factors.
        Then $P = t_1 d_1 = d_n t_n$, and $t_i d_i = d_{i-1} t_{i-1}$
        when $i>1$.  

        Let $k_1$, $k_2$, $\ldots$, $k_n$ be arbitrary integers.  Then 
        the series
        \[
                [t_1, d_1]_{k_1} + q^{-k_1} [t_2, d_2]_{k_2} + 
                        q^{-k_1-k_2} [t_3, d_3]_{k_3} + \cdots +
                        q^{-k_1 - \cdots - k_{n-1}}[t_n, d_n]_{k_n}
        \]
        is a telescoping series which reduces to 
        \[      t_1 d_1 - q^{- \sum k_i} d_n t_n 
                        = (1-q^{-\sum k_i}) P.
        \]

        By the induction hypothesis, there are $T_i$ homogeneous in 
        the $\partial^{\beta^{a}}$'s of degree $n$ such 
        that $[T_i, x] = t_i \sigma_b$.  
        For each $i = 1, \ldots, n$, put
        $k_i = -n a_i$.
        The homogeneity of the $T_i$ ensures that for each $i$, 
        \[      T_i [\partial^{\beta^{a_i}}, x] = T_i \sigma_{a_i}
                                = q^{na_i} \sigma_{a_i} T_i
                                = [\Si_{-na_i}(\partial^{\beta^{a_i}}), x] T_i, 
        \]
        Thus we get 
        \[      T_i [d_i, x] - [\Si_{k_i}(d_i),x] T_i = 0.
        \]
        We can use this and (\ref{tim:2}) to get 
        \[      [[T_i, d_i]_{k_i}, x] = [t_i\sigma_b, d_i]_{k_i}
        \]
        
        Thus, applying $[\cdot, x]$ to 
        \[
                \tilde{Q}=[T_1, d_1]_{k_1} + q^{-k_1} [T_2, d_2]_{k_2} + 
                        q^{-k_1 - k_2} [T_3, d_3]_{k_3} + \cdots +
                        q^{-(k_1 + \cdots + k_{n-1})}[T_n, d_n]_{k_n}
        \]
        will yield  
        \begin{align}\label{tim:4}
                        [t_1 \sigma_b, d_1]_{k_1} 
                        + & q^{-k_1} [t_2 \sigma_b, d_2]_{k_2} + 
                        q^{-(k_1+k_2)} [t_3 \sigma_b, d_3]_{k_3} + \\ 
                        & \cdots +
                        q^{-(k_1 + \cdots + k_{n-1})}[t_n \sigma_b, d_n]_{k_n}.
                        \notag
        \end{align}
        We use the identity
        $[t_i \sigma_b, d_i]_{k_i} = q^{-b}[t_i, d_i]_{k_i-b} \sigma_b$
        on (\ref{tim:4}) to get
        \begin{align*}
                        q^{-b}\bigg( [t_1, d_1]_{k_1 - b}
                        + & q^{-k_1}[t_2, d_2]_{k_2-b} + 
                        q^{-(k_1+k_2)}[t_3, d_3]_{k_3-b} + \\ 
                         & \cdots + q^{-(k_1 + \cdots + k_{n-1})}
                                [t_n, d_n]_{k_n-b} \bigg) \sigma_b.
        \end{align*}
        which can be reduced to 
        \[
                q^b(1-q^{-\sum(k_i - b)})P 
                        \sigma_b.
        \]
        We chose the $k_i$ so that $\sum k_i = - n \sum a_i$.  Thus 
        $-\sum(k_i - b) = n (b + \sum a_i)$.
        By our assumption on $b$, either all $a_i =0$ or 
        $c =q^b(1-q^{n(b+ \sum a_i)})$ is nonzero.

        If some $a_i \neq 0$, then
        $[c^{-1}\tilde{Q}, x] = P$ as required.
        However, if all $a_i = 0$, then $P = \partial^n$.  
        In this case,
        $Q = \frac{1}{n+1} \partial^{n+1}$ is the required operator.
\end{proof}

        Now we are now ready to prove our main theorem.

\begin{proof}[Proof of Theorem~\ref{T:qdops}]
        Let $A^n$ be the $\K$-module generated by 
        \[
                \{ \partial^{\beta^I}\sigma_a \mid 
                \text{ $|I| \leq n$ and $a \in \Z$ } \}.
        \]
        Since $\sigma_a$ commutes with $\partial^{\beta^I}$ up to 
        multiplication by a scalar in $\K$, we have 
        \[
        D^0_q A^n D^0_q = R A^n R
        \] for any $n$.  
        We will show by induction on $n$ that $R A^n R = D^n_q$. 
        Since $A^0 = \K \langle \sigma_a \mid a \in \Z \rangle$, 
        we have $R A^0 R = D_q^0$.  This proves the base case.

        Now suppose that $R A^{n-1} R = D^{n-1}_q$.  By 
        Corollary~\ref{C:2} above, we know that 
        $A^n \subset D^n_q$.
        Thus, we must show $D^n_q \subseteq R A^n R$.
        Since  $D^n_q = D^0_q \CZ_q^{'n} D^0_q$, and
        $R A^n R = D^0_q A^n D^0_q $, it suffices to show 
        $\CZ_q^{'n} \subseteq R A^n R$.

        Let $\varphi \in Z_q^{'n}$.  Then $[\varphi, x] \in R A^{n-1} R$
        by the induction hypothesis.
        Thus $[\varphi, x]= \sum r_j P_j \sigma_{b_j} s_j$ where 
        $r_j, s_j \in R$, $P_j = \partial^{\beta^{I_j}}$ with
        $|I_j| \leq n-1$, and $b_j \in \Z$. 
        By Lemma~\ref{L:intmono}, we can find $Q_j \in A^n$ such 
        that $[Q_j, x] = P_j \sigma_{b_j}$.  Put $Q = \sum r_j Q_j s_j$.
        Then $Q$ is in $R A^n R$.
        Since $r_j$ and $s_j$ commute with $x$, we have
        $[Q,x]= \sum r_j P_j \sigma_{b_j} s_j$.  Hence $[\varphi - Q , x]=0$.
        That is, $\varphi - Q = \eta \in D^0_q \subset R A^n R$.  Hence, 
        $\varphi = Q + \eta$, and so $\varphi \in R A^n R$.

        We have shown that $A = \bigcup A^n$ generates $D_q$ over
        $R$.  
        Moreover, it is clear from the construction of $A^n$ that
        $A$ is generated as a $\K$-algebra by 
        $\{ \partial^{\beta^a} \mid a \in \Z \}$
        and $\{ \sigma_a \mid a \in \Z \}$.
        Furthermore, by Note~\ref{N:nto1}, 
        each $\partial^{\beta^a}$ is in the $\K$-span of
        $\{ \sigma_a \mid a \in \Z \} 
        \bigcup \{\partial , \partial^{\beta} \}$.
        Hence, it is enough to show that each $\sigma_a$ is generated
        over $R$ by 
        $\partial^{\beta^{-1}}$, $\partial$, and $\partial^{\beta}$.

        To this end, note that 
        $\sigma_{-1} = \partial^{\beta^{-1}} x - x \partial^{\beta^{-1}}$, 
        $\sigma_0 = 1 = \partial x - x \partial$, and 
        $\sigma_{1} = \partial^{\beta} x - x \partial^{\beta}$.
        For any positive integer $a$, $\sigma_{-a} = (\sigma_{-1})^a$
        and $\sigma_a = (\sigma_1)^a$.
        This proves the theorem.
\end{proof}
\begin{remark}\label{R:betasigma}
In remark 3.1.9 of \cite{LR1}, Lunts and Rosenberg show that
$D_q \supset D_{\beta}\sigma ({\Gamma})$, 
where $D_{\beta}$ denotes the ring of 
$\beta$ differential operators.  In our case, $D_{\beta}$ is the 
$\Bbbk$-algebra
generated by $x,\partial ^{\beta}$ with relations 
$[\partial ^{\beta} , x]_1 = 1$ (for details, see \emph{Construction of
Skew Polynomial Rings}, page 7 of \cite{GW}).
Our work shows that $D_q \neq D_{\beta}\sigma ({\Gamma})$. Indeed, the usual
derivation $\partial \notin D_{\beta}\sigma ({\Gamma})$. 
By degree considerations,
if $\partial \in D_{\beta}\sigma ({\Gamma})$, then
$\partial = f \partial ^{\beta}$ where $f$ is of degree 0.
Hence, $f = \sum _{-n \leq i \leq m} c_i \sigma ^i$. 
As $f (x^r) = \frac{r}{(1 + q + \cdots + q^{r-1})} 
x^r$ for $r\geq 1$, we can choose $r >> 0$
such that $\sum_{-n\leq i \leq m} c_i q^{ri} (1 + q + \cdots + q^{r-1}) = r$
implies that $c_{-n} ,c_m = 0$.
\end{remark}
\begin{theorem}\label{T:Dqsimple}
The ring $D_q$ is a simple ring.
\end{theorem}
\begin{proof}
Let $\mathcal{I}$ be an ideal in $D_q$.  Let $0 \neq f\in \mathcal{I}$.
Then $f $ can be written as
\[
f = \sum _{\{ a,I\mid |I| \gneq 0\} } \sigma (a) p_I(x) \partial ^{\beta ^I}
        + \sum _{n\in \mathbb{Z}} \sigma (n) p_n(x),
\]
where $p_I, p_n$ are polynomials in $x$, $a \in \mathbb{Z}$ 
and the multi indices $I$ have entries
in $\{ 0,1 \}$.
We induct on $d = \mathit{max} \{ |I| \mid p_I \neq 0 \}$, 
to claim that the ideal 
containing $f$ should contain 1.
Assume that $d= 0$.
Then 
\[
f = \textit{lower degree in } x + (c_1\sigma (a_1) + \cdots + c_r 
        \sigma (a_r))x^m,
\]
where $c_i \in \Bbbk$ and $a_1\lneq a_2 \lneq \cdots \lneq a_r$.
Now $[\sigma (-a_r) f ,x]$ have fewer monomials than $f$.
Continuing thus, we can assume that $f = c x^b$ for some $c \in \Bbbk$.
Since $[\partial ,x^b] = bx^{b-1}$ the base case is proved.

Assume that the claim has been proved for all positive integers less than
$d$. We can write $f$ as
\[
f = \textit{lower lengths in }I + \sigma (a_{1}) p_{I_1}(x) 
        \partial ^{\beta ^{I_1}} + \cdots +
        \sigma (a_{r})p_{I_r}(x)\partial ^{\beta ^{I_r}},
\]
where
$a_{s} \in \mathbb{Z}, a_{1} \lneq a_{2} \lneq \cdots \lneq a_{r}$
and $p_{I_s}(x)$ are polynomials in $x$ where $I_s$ are multi indices of 
length $d$, and these multi indices can repeat in the above sum.
This implies that,
\[
[\sigma (-a_{r}) f , x] = \textit{lower lengths in }I +
        \sigma (b_{1}) g_{I_1}(x)\partial ^{\beta ^{I_1}} + \cdots +
        \sigma (b_{{r-1}}) g_{I_{r-1}}(x)\partial ^{\beta ^{I_{r-1}1}}
\]
where $b_{s} \in \mathbb{Z}, g_{I_l} \in \Bbbk [x]$ and 
$b_{1} \lneq b_{2} \lneq \cdots \lneq b_{{r-1}}$.
This shows that $[\sigma (-a_{r}) f , x]$ has fewer monomials which
correspond to multi indices of length $d$. This completes the theorem.
\end{proof}


\section{An intrinsic description of $D_q$}\label{sect:intrinsic}

Although we know that the ring $D_q$ is generated over $R$
by $\partial^{\beta^{-1}}$, $\partial$, and  $\partial^{\beta}$, 
the relations among these are not necessarily apparent.  
In this section, we will give a description of $D_q$ 
in terms of generators and relations.

We begin with a study of the grading on $D_q$.
First, to simplify notation, let $\tau = x \partial$ and $\sigma = \sigma_1$.
Then $\sigma^{-1} = \sigma_{-1}$.
\begin{lemma}\label{L:D0domain}
        The ring $(D_q)_0$ of homogeneous $q$-differential 
        operators of degree 0 is $\K[\sigma, \tau, \sigma^{-1}]$. 
\end{lemma}
\begin{proof}
        It is clear that $\sigma$, $\sigma^{-1}$, and $\tau$ are 
        homogeneous of degree 0.  We must show that they generate
        all of $(D_q)_0$, that they commute, and that they have
        no other relations among them.

        Suppose that $\phi \in (D_q)_0$.  Then by Theorem~\ref{T:qdops},
        $\phi$ can be written as $\phi = \sum P_i$ where each $P_i$ is 
        a monomial in $x$, $\partial^{\beta^{-1}}$, $\partial$, and  
        $\partial^{\beta}$.  
        Moreover, each $P_i$ has degree $0$, so each $P_i$ has exactly 
        one $x$ for each $\partial^{\beta^{-1}}$, $\partial$, and  
        $\partial^{\beta}$.
        Since $\partial^{\beta^{a}} x = q^a x \partial^{\beta^a} + 1$ for 
        $a = -1, 0, 1$, we can rewrite $\phi$ so that each $P_i$ has 
        the form 
        $P_i = (x \partial^{\beta^{a_1}}) \cdots (x \partial^{\beta^{a_k}})$.
        If $a_j \neq 0$, then 
        $(\sigma_{a_j} - 1)/ (q^{a_j} - 1) = x \partial^{\beta^{a_j}}$.
        Since $\tau = x \partial^{\beta^0}$, we have that each $P_i$ is in 
        the span of $\sigma$, $\sigma^{-1}$, and $\tau$.  Hence 
        $(D_q)_0$ has the required generators.

        Since $\sigma$ and $\sigma^{-1}$ are multiplicative inverses, 
        they commute.  Comparing the values of $\tau \sigma$ and 
        $\sigma \tau$ on $x^m$ for any $m$ shows that $\sigma$ and $\tau$ 
        also commute. 
        All we need to show now is that there are no remaining relations.

        Suppose that we have a relation
        $\sum a_{ij} \sigma^i \tau^j = 0$ for some $a_{ij} \in \K$.
        Then we can multiply this expression by $q$ and 
        $\sigma$ an appropriate number of times to ensure that all 
        all powers of $q$ in the $a_{ij}$ and all $i$ are nonnegative.
        Then for every $m$, 
        $0 = \sum a_{ij} \sigma^i \tau^j (x^m) 
           = (\sum a_{ij}q^{im} m^j) x^m$.
        Fix an $m >> 1$ such that for every $i$ and $j$, all the powers 
        of $q$ appearing in $a_{ij}$ are less than all the powers of $q$ 
        appearing in $a_{(i+1)j}q^m$.  Then the polynomial 
        $\sum a_{ij}q^{im} m^j$ is a polynomial in $q$ with no 
        terms canceling.  Since it is zero, all its coefficients
        are zero.  But if $\lambda$ is the coefficient of $q^k$ in
        $\alpha_{ij}$, then $m^j \lambda $ is the coefficient of 
        $q^{im + k}$ in  $\sum a_{ij}q^{im} m^j$.
        Hence $\lambda = 0$.  It follows that $a_{ij} =0$ 
        completing the proof.
\end{proof}

\begin{corollary}\label{C:Dqdomain}
        The ring $D_q$ is a domain.
\end{corollary}
\begin{proof}
        First we show that $x$ cannot be a zero divisor.  
        It is obvious that
        $x \cdot \varphi \neq 0$ if $\varphi \neq 0$.
        Suppose $\varphi \cdot x = 0$, then 
        $\varphi (x^n) = 0$ for all
        positive integers $n$.  Also note that 
        $[\varphi ,x]_a = -q^a x\varphi$.  Hence $\varphi$ does not
        satisfy Lemma \ref{L:comm}.

        Suppose $\phi \psi =0$.  Let $\phi_a$ and $\psi_b$ be the 
        highest degree parts of $\phi$ and $\psi$ of degrees $a$ and
        $b$ respectively.  
        Then $\phi_a \psi_b = 0$. 
        
        If $b$ is positive, then $(\phi _a x^b) \psi ^{\prime}=0$
        where $\psi ^{\prime}$ is of degree 0.
        Note also that $\phi _a x^b$ is a non-zero homomorphism.
        If $b$ is negative, then $\phi _a (\psi _b x^{-b}) = 0$
        and $\psi _b x^{-b} \neq 0$ is of degree 0.
        Hence, we can assume that $b=0$.  

        Now, if $a$ is positive, then
        $x^a \phi ^{\prime} \psi _b = 0$ implies 
        $\phi ^{\prime} \psi _b = 0$, where $\phi ^{\prime}$ is of degree 0.
        If $a$ is negative then 
        $(x^{-a} \phi _a) \psi _b = 0$ where $x^{-a} \phi _a$ is of
        degree 0.  Hence, we can assume that $a=0$ also.
        
        Now the corollary follows from the lemma \ref{L:D0domain}.
\end{proof}
Now we can describe the ring $D_q$ intrinsically.
\begin{theorem}\label{T:intrinsic}
        The ring $D_q$ of $q$-differential operators on $R$ is 
        the $\K$-algebra generated by $x$, $\partial^{\beta^{-1}}$, 
        $\partial^{\beta^0}$, and  
        $\partial^{\beta^1}$ subject to the relations
        \begin{align*}
                \partial^{\beta^{a}} x - q^{a} x \partial^{\beta^{a}}
                        &= 1, &
                \partial^{\beta^{a}} x \partial^{\beta^{b}} 
                        &= \partial^{\beta^{b}} x  \partial^{\beta^{a}}, &
                        &\text{and} &
                \partial^{\beta^{-1}} - q  \partial^{\beta} 
                        &= (1 - q)\partial^{\beta^{-1}} x  \partial^{\beta}.
        \end{align*}
\end{theorem}
\begin{proof}
        Let $F$ be the $\K$-algebra generated by symbols $x$, 
        $\partial^{\beta^{-1}}$, $\partial^{\beta^{0}}$, 
        and $\partial^{\beta^{1}}$ subject to the relations 
        \begin{equation}\label{tim:sect2.1}
                \partial^{\beta^{a}} x - q^{a} x \partial^{\beta^{a}} = 1
                \text{\ for \ } a = -1, 0, 1.
        \end{equation}
        Giving each $\partial^{\beta^{a}}$ degree -1 and giving $x$ degree 1
        makes $F$ into a graded $\K$-algebra.  In fact, the natural 
        quotient map $\pi: F \to D_q$ preserves this grading.
        Let $I$ be the kernel of $\pi$.  Then $I$ is generated by 
        homogeneous elements.

        Suppose that $\theta \in I$ is a homogeneous generator of $I$
        of degree $n$.  If $n > 0$, then $\theta \partial^n$ is a 
        homogeneous element of $I$ of degree 0.
        If $n < 0$, then $x^{-n} \theta$ is likewise a homogeneous 
        element of $I$ of degree $0$.  Since neither $x$ nor $\partial$
        are zero divisors in $F$, neither $x^{-n} \theta$ nor 
        $\theta \partial^n$ are zero.  Hence every homogeneous generator
        of $I$ is a factor of some degree $0$ element of $I$.

        Now let us restrict our attention to degree $0$.  The map 
        $\pi$ takes the ring $F_0$ of degree 0 elements of $F$
        to the ring $(D_q)_0$ with 
        kernel $I_0$.  Using (\ref{tim:sect2.1}), we can write every 
        element $F_0$ as a polynomial in $x\partial^{\beta^{-1}}$,
        $x\partial^{\beta^{0}}$, and $x\partial^{\beta^{1}}$.
        Since $(D_q)_0$ is commutative, the commutators of any two
        elements of $F_0$ must be in $I_0$.  Hence, 
        \[ 
                x \partial^{\beta^a} x \partial^{\beta^b} - 
                x \partial^{\beta^b} x \partial^{\beta^a} \in I_0
        \]
         for any $a$ and $b$.  

        The only remaining identity in $(D_q)_0$ comes from the 
        fact that $\sigma_1$ and $\sigma_{-1}$ are multiplicative 
        inverses.  Since $\sigma_a = 1 + (q^a - 1)x \partial^{\beta^a}$
        for any $a$, we have
        \begin{align*}
            1 &= \sigma_1 \sigma_{-1} \\
                &= (1 + (q - 1) x \partial^{\beta^{1}}) 
                  (1 + (q^{-1} - 1) x \partial^{\beta^{-1}})\\
                &= 1 + (q-1) x 
                  (\partial^{\beta^1} - q^{-1} \partial^{\beta^{-1}}
                + (q^{-1} - 1) \partial^{\beta^{1}} x \partial^{\beta^{-1}})
        \end{align*}
        Hence, 
        \[ 
                (\partial^{\beta^1} - q^{-1} \partial^{\beta^{-1}}
                + (q^{-1} - 1) \partial^{\beta^{1}} x \partial^{\beta^{-1}}) 
                \in I_0.
        \]
        
        Let us recapitulate:
        \noindent
        In $(D_q)_0$, we have only the standard relations
        \[ \partial^{\beta^{a}} x - q^{a} x \partial^{\beta^{a}} = 1, \]
        the commutator relations
        \[ x \partial^{\beta^{a}}  x \partial^{\beta^{b}} 
                =  x \partial^{\beta^{b}}  x \partial^{\beta^{a}}, \]
        and the special relation 
        \[ x(\partial^{\beta^1} - q^{-1} \partial^{\beta^{-1}}
          =(1 - q^{-1}) \partial^{\beta^{1}} x \partial^{\beta^{-1}}). \]  
        Since all relations in $D_q$ come from elements in $I$, and 
        all generators of $I$ are factors of elements in $I_0$, 
        the only other possible relations in $D_q$  are 
        $\partial^{\beta^{a}}  x \partial^{\beta^{b}} 
                = \partial^{\beta^{b}}  x \partial^{\beta^{a}}$
        and
        $\partial^{\beta^1} - q^{-1} \partial^{\beta^{-1}}
           = (1- q^{-1}) \partial^{\beta^{1}} x \partial^{\beta^{-1}}$. 
        Inspecting these last two expressions on $x^m$ shows that 
        they do indeed hold. 
\end{proof}
The following formulae are immediate, and hence we do not provide any
proofs:\hfill
\begin{enumerate}
        \item   $
                x [\partial , \partial ^{\beta}]_1 
                = \partial -\partial ^{\beta}$;
                $[\partial , \partial ^{\beta}]_1 x= \partial - q 
                \partial ^{\beta}$.
        \item   $(\tau + k) \partial ^{\beta ^a} = \partial ^{\beta ^a} 
                        (\tau + k -1)$.
        \item   $(\tau + 1) \partial ^{\beta} = (\frac{q\sigma -1}{q-1})
                        \partial$.

                This can be generalized to multi-indices 
                $I = (i_1,i_2,\cdots ,i_n)$ 
                (using the notations as in \ref{E:multiindex})
                with entries
                in $\{ 0,1 \}$ as
                \[
                \left( \prod _{ \{j|i_j = 1 \}} (\tau +j) \right)
                \partial ^{\beta ^I} =
                \left( \prod _{\{j|i_j = 1\} } 
                \frac{q^j\sigma - 1}{q-1} \right) \partial ^{|I|}
                \]
\end{enumerate}
\begin{remark}
We have not been successful in determining whether $D_q$ is left noetherian 
or not. 
\end{remark}
\section{Generalization to several variables}\label{sect:nvariables}
Let $q_1,q_2,\cdots ,q_n$ be transcendental elements over
$\mathbb{Q}$ and let $\K $ be a field containing $\mathbb{Q} (q_1,q_2,\cdots,
q_n)$.   
Let $R = \K[x_1,x_2,\cdots ,x_n]$, and $D^m_q$ denote $D^m_q(R)$ (respectively
$D_q$ denote $D_q(R)$).  The ring $R$ has a $\mathbb{Z}^n$-grading given by
\[
\mathit{deg}(x_1^{a_1}x_2^{a_2} \cdots x_n^{a_n}) = (a_1,a_2,\cdots ,a_n).
\]
Let $\beta : \mathbb{Z}^n \times \mathbb{Z}^n \to \K^{\times}$ be a bicharacter
defined by
\[
\beta (\mathbf{a} ,\mathbf{b}) = q_1^{a_1b_1}q_2^{a_2b_2}\cdots q_n^{a_nb_n}.
\]
For $\mathbf{a} = (a_1,a_2,\cdots ,a_n) \in \mathbb{Z}^n$ let
 $\mathbf{x^a} := x_1^{a_1}x_2^{a_2} \cdots x_n^{a_n}$.
For each $i,1\leq i \leq n$, define $\partial _i ^{\beta ^k}$ as
\[
\partial _i ^{\beta ^k} (\mathbf{x^a}) = \frac{(q_i^{ka_i} -1)}{(q_i - 1)}
                        \mathbf{x}^{\mathbf{a}-(0,0,\cdots ,1_i,0,\cdots ,0)}.
\]
Note that 
\begin{align*}
[\partial _i ^{\beta ^k} , x_j] & = 0 \textit{ for } i \neq j,\\
[\partial _i^{\beta ^k},\partial _j^{\beta ^m}] &=0  \textit{ for } i \neq j,\\
[\partial _i ^{\beta ^k} ,\sigma _{\mathbf{a}}] &= 0 \textit{  when }a_i = 0.
\end{align*}
The notes \ref{N:nto1}, lemma \ref{L:dinD1} and corollary \ref{C:2}
follow verbatim.
\begin{theorem}\label{T:nvariables}
The ring $D_q$ of $q$-differential operators on $R$ is generated as
a $\K$-algebra by $R$ and the set $\{ \partial _i ^{\beta ^{-1}},\partial _i,
\partial _i ^{\beta} \}_{1\leq i \leq n}$.
\end{theorem}
\begin{proof}
Let $A$ denote the 
$\K$-algebra generated 
by $R$ and the set $\{ \partial _i ^{\beta ^{-1}},\partial _i,
\partial _i ^{\beta} \}_{1\leq i \leq n}$.
Given  $f_i\in R$ and $P_i$ in terms of monomials consisting of 
$\{ \partial _j ^{\beta ^k} | k,j \in \mathbb{Z} \}$,
let
$Q$ be such that
$[Q,x_i] = P_i \sigma _{\mathbf{a}_i}$.
We show that $Q\in A$.
We can think of $[. , x_i] $ as $\partial _{y_i}$, where
\begin{align*}
\partial _{y_i} (\partial _j) &= \delta _{i,j},\\
\partial _{y_i} (\partial _j ^{\beta}) &=  
                \delta _{i,j} \sigma _{(0,0,\cdots ,1_i,0,\cdots 0)},\\
\partial _{y_i} (x_j) &=0, \textit{  and }\\
\partial _{y_i} (\sigma _{\mathbf{a}}) &= 
        (q^{a_i} - 1) x_i \sigma _{\mathbf{a}}.
\end{align*}
By lemma \ref{L:intmono},
for each $i$, we can find 
a '$y_i$-integral' $Q_i$ of $f_iP_i \sigma _{\mathbf{a}_i}$ in $A$.
Since
\[
[[Q,x_i],x_j] = [[Q,x_j],x_i],
\]
we have $\partial _{y_i} \partial _{y_j} = \partial _{y_j} \partial _{y_i}$.
Thus,
we can find an $F \in A$ 
such that
$[Q-F, x_i] = 0$ for all $i$.  That is, $Q-F \in D^0_q$ which is contained in
$A$.
Hence the theorem.
\end{proof}
\begin{remark}
The ring $D_q$ is simple and a domain.  
\end{remark}
\section{Relationship with the Quantum group on $sl_2$}\label{S:U_q(sl2)}
\subsection{The ring $\Gamma_q(\mathbb{P}^1)$}\label{S:globalsections}
Fix $\beta : \mathbb{Z} \times \mathbb{Z} \to \K ^{\times}$
be fixed as $\beta (n,m) = q^{nm}$. For this $\beta$, we define the following
rings.
Let $D_q $ (respectively $L_q$) 
denote the ring of
$q$-differential operators on $\K [x]$ (respectively
$\K [y]$)
with the set-up 
as in section \ref{S:1var} (respectively, degree of y = -1).
$L_q$ is a \K -algebra generated 
by $\partial _y ,\partial ^{\beta}_y,\partial ^{\beta ^{-1}}_y $
where  $\partial _y ^{\beta}$ is
a left $\sigma _y$-derivation for $\sigma _y (y) = \frac{1}{q}y$ and
$\partial _y^{\beta ^{-1}}$ is a left $\sigma _y^{-1}$-derivation.
Hence we have the following formulae:
\begin{align*}
\partial _y (y^n) &= n y^{n-1},\\
\partial _y ^{\beta} (y^n) &= \left( \frac{1- \frac{1}{q^n}}{1-\frac{1}{q}}
                        \right) y^{n-1},\\
\partial _y^{\beta ^{-1}}(y^n) &= \left(\frac{1-q^n}{1-q} \right) y^{n-1}.
\end{align*}

Similarly, let $M_q$ denote the
ring of $q$-differential operators on $\K [x,x^{-1}]$ where $\K [x,x^{-1}]$
is $\mathbb{Z}$-graded as $\mathit{deg}(x) = 1$ and
$\mathit{deg}(x^{-1}) = -1$.
By theorem 3.2.2 of \cite{LR1}, 
there are canonical ring homomorphisms (injective) $D_q \to M_q$ and
$L_q \to M_q$. Specifically, the maps are
\begin{align*}
D_q &\to M_q & L_q &\to M_q\\
\sigma ^{\pm} (y^n) &= q^{\mp n}y^n & \sigma _y^{\pm} (x^n) &= q^{\pm n} x^n,\\
\partial (y)&= -y^2 & \partial _y (x) &= -x^2,\\
\partial ^{\beta}(y) &= -\frac{1}{q} y^2 & \partial ^{\beta}_y (x) &= -qx^2,\\
\partial ^{\beta ^{-1}} (y) &= -qy^2 & \partial ^{\beta ^{-1}}_y (x) &=
                                                -\frac{1}{q}(x^2),
\end{align*}
and extend the respective derivations (or $\beta$-derivations) to the entire
ring.
We define $\Gamma _q(\mathbb{P}^1)$ as 
\begin{align*}
\Gamma _q (\mathbb{P}^1) := \mathit{Ker} \hookrightarrow D_q \bigoplus L_q
                         &\to M_q;\\
(\varphi _1 ,\varphi _2) &\mapsto \varphi _1 - \varphi _2.
\end{align*}
For simplicity, let $\Gamma _q := \Gamma _q (\mathbb{P}^1)$.
Note that $\Gamma _q$ is a ring because the homomorphisms
$D_q \to M_q$ and $L_q \to M_q$ preserve multiplication.

\begin{lemma}
$\Gamma _q $ is generated 
 over $\K$
by the set 
\begin{align*}
\{ (\partial  ,-y^2\partial _y),
(-x^2\partial ,\partial _y), 
(\partial ^{\beta}, -\frac{1}{q}y^2\partial _y ^{\beta}),
(-qx^2\partial ^{\beta},\partial _y ^{\beta}),\\
(\partial ^{\beta ^{-1}}, -qx^2\partial _y ^{\beta ^{-1}}),
(-\frac{1}{q}x^2\partial ^{\beta ^{-1}},\partial _y ^{\beta ^{-1}})\} .
\end{align*}
\end{lemma}
\begin{proof}
Let the algebra generated over $\K$ 
by the set mentioned in the statement of the
lemma be $G$.

Clearly, members of the set
$\{ \partial ,\partial ^{\beta},\partial ^{\beta ^{-1}},
        x^2\partial ,x^2 \partial ^{\beta},x^2\partial ^{\beta ^{-1}} \}$
map $\K [x^{-1}] $ to itself.
Note that 
\begin{align*}
\sigma &= \left( \frac{q-1}{q+1}\right) 
                [\partial ^{\beta} , x^2\partial ^{\beta}]_2 +1,\\
\sigma ^{-1} &= \left( \frac{\frac{1}{q}-1}{\frac{1}{q}+1}\right)
                [\partial ^{\beta ^{-1}},
                        x^2\partial ^{\beta ^{-1}}]_{-2} +1,\\
\tau = x\partial &= \frac{1}{2}[\partial , x^2\partial ].
\end{align*}
Thus, elements of $(D_q)_0$ (degree 0 endomorphisms in $D_q$) 
map $\K [x^{-1}]$ to itself.
If $(\varphi _1,\varphi _2) \in \Gamma _q$,
then $\varphi _1 (\K [y]) \subset \K [y]$ and
$\varphi _2 (\K [x]) \subset \K [x]$ and 
$\mathit{deg} (\varphi _1) =\mathit{deg}(\varphi _2)$.
If $n = \mathit{deg}(\varphi _1) \leq 0$, then
$\varphi _1 = \sum _{|I| = -n} f_I \partial ^{\beta ^I}$,
where $f_I \in (D_q)_0$.
Therefore $\varphi _1$ is generated by $\{ \partial ,\partial ^{\beta} ,
        \partial ^{\beta ^{-1}} \}$.
Thus,  $(\varphi _1, \varphi _2 ) \in G$.
Similarly, if $n = \mathit{deg}(\varphi _2) \geq 0$,
then $\varphi _2  = \sum _{|I| =n}g_I \partial _y ^{\beta ^I}$ where
$g_I \in (D_q (\K [y]))_0$.  Again $(\varphi _1 ,\varphi _2) \in G$.
Hence the lemma.
\end{proof}
\subsection{Quantum group on $sl_2$}\label{S:UqtoDq}
Let $U_q$ denote the Quantum group corresponding to the Lie algebra
$sl_2(\K )$.  That is, $U_q$ is a $\K$-algebra generated by $E,F,K,K^{-1}$,
with relations given by (for details see\cite{CP} or \cite{J})
\begin{align*}
KK^{-1} = &1 = K^{-1}K, \\
KEK^{-1} &= q^2E,\\
KFK^{-1} &= q^{-2}F,\\
EF - FE &= \frac{K - K^{-1}}{q-q^{-1}}.
\end{align*}
Recall that $q$ is a transcendental element over $\mathbb{Q}$ throughout 
this paper, and $\K$  
contains $\mathbb{Q}(q)$.
The ring $U_q$ is a Hopf-algebra and for the purposes of this paper, we will 
give the comultiplication map $\Delta$:
\begin{align*}
\Delta (E) &= E\otimes 1 + K \otimes E,\\
\Delta (F) &= 1 \otimes F + F \otimes K^{-1},\\
\Delta (K) &= K \otimes K.
\end{align*}
Let the ring corresponding to the \textit{quantum plane} be $S$.
That is,
\[
S = \K < u ,v > \diagup uv = qvu.
\]
There is an action of $U_q$ on this ring given by
\begin{align*}
K(1) &=1   & E(1) &= F(1) = 0,\\
K(u) &= qu & K(v) &= \frac{1}{q} v,\\
E(u) &= 0  & E(v) &= u,\\
F(u) &= v  & F(v) &= 0,
\end{align*}
and extend the action on $S$ via $\Delta$.  If we consider the Ore set
$\{ v^n \}_{n\geq 0}$
in $S$, then this action of $U_q$ extends (via $\Delta$) to the 
Ore-localization.
For example, the extension of $E$ to the Ore-localization of $S$ is as
follows:
\begin{align*}
0 = E(1) &= E(v \frac{1}{v}) \\
        &= E(v) \frac{1}{v} + K(v) E(\frac{1}{v})\\
        &= u\frac{1}{v} +\frac{1}{q} v E(\frac{1}{v}).
\end{align*}
Hence, $E(\frac{1}{v}) = -q \frac{1}{v} u \frac{1}{v}$. 
This extended action of $U_q$ keeps the
polynomial ring
$\K [u\frac{1}{v}]$ invariant.
We let $x : = u\frac{1}{v}$. Then we have
a homomorphism of \K -algebras
\[
\alpha : U_q \to D_q ,
\]
given by
\begin{align*}
\alpha (F) &= q^{-1}\sigma ^{-2}\partial ^{\beta ^2} = q^{-1}\sigma ^{-2}
                \left( \frac{1+q\sigma}{1+q} \right) \partial ^{\beta},\\
\alpha (E) &= -q^2 x^2 \partial ^{\beta ^2} = -q^2 \left( 
                \frac{1+q^{-1}\sigma}{1+q}\right)
                 x^2 \partial ^{\beta},\\
\alpha (K) &= \sigma ^2,\\
\alpha (K^{-1}) &= \sigma ^{-2}.
\end{align*}
Note that $\alpha (U_q) \subset D_{\beta} \sigma (\mathbb{Z})$ (referred to in
remark \ref{R:betasigma}) and hence $\alpha$ is not surjective.
Similarly, there is an algebra homomorphism
\[
\gamma : U_q \to L_q,
\]
given by
\begin{align*}
\gamma (F) &= -\frac{1}{q^2} 
                \sigma _y^{-2} y^2\partial ^{\beta ^2}_y= - 
                                        \frac{1}{q^2}\sigma _y^{-2}
                        \left( \frac{1+q\sigma _y}{1+q}\right) 
                        y^2\partial _y^{\beta},\\
\gamma (E) &= q \partial ^{\beta ^2}_y = q\left( 
                \frac{1+q^{-1}\sigma _y}{1+q}\right)\partial _y^{\beta},\\
\gamma (K) &= \sigma _y ^{2},\\
\gamma (K^{-1}) &= \sigma _y ^{-2}.
\end{align*}
Again, $\gamma$ is not a surjection.
The above two homomorphisms give a homomorphism
\begin{align*}
\eta :U_q &\to \Gamma _q;\\
u &\mapsto (\alpha (u),\gamma (u)).
\end{align*}
Since $(\partial , -y^2\partial _y) \notin \eta (U_q)$, we have the following
\begin{proposition}
$\eta$ does not give a surjection of $U_q$ to $\Gamma _q$.
\end{proposition}
But we have a surjection by considering inverse limits. This is shown in the
following two  subsections.
\subsection{Inverse limits of $\Gamma _q$}
Let \CA $=\mathbb{Q}[q,q^{-1}]_{(q-1)}$.
Let $D_{q,\CA}(\CA /(q-1)^n[x])$ 
denote the ring of $\CA /(q-1)^n$ -linear $q$-differential operators
on $\CA /(q-1)^n [x]$.
Let $D_t(\mathbb{Q}[x,t] /t^n)$ denote the ring of $\mathbb{Q} [t]/t^n$-linear
usual differential differential operators on $\mathbb{Q}[x,t]/t^n$.
That is, $D_t (\mathbb{Q}[x,t]/t^n) = \mathbb{Q}[t]<x,\partial > \diagup
                                     t^n,[\partial ,x] =1 $.
\begin{lemma}
The rings $D_{q,\CA}(\CA /(q-1)^n[x])$ and 
 $D_t(\mathbb{Q}[x,t] /t^n)$ are isomorphic as $\mathbb{Q}$-algebras.
\end{lemma}
\begin{proof}
First we note that
\begin{align*}
(\CA /(q-1)^n)[x] &\cong \mathbb{Q}[x,t]/t^n;\\
(q-1) &\mapsto t,
\end{align*}
as $\mathbb{Q}$-algebras.
So, it suffices to show that if $d \in D_{q,\CA}(\CA /(q-1)^n[x])$ 
then 
\[
[\cdots [[d,x],x],\cdots x] = 0
\]
 for some finite number of commutators.
Suppose $\varphi \in D^0_{q,\CA}(\CA/(q-1)^n [x])$ is such that
$[\varphi , x]_m = 0$ for some $m$.  Then for any $c\in \CA /(q-1)^n [x]$,
we have $[c\varphi ,x] = (q^m-1)xc\varphi$. Thus, the $n$-commutators
$[\cdots [[c\varphi ,x],x]\cdots ]=0$.  

Now suppose that $[\varphi ,x]_m = d \in D^{l-1}_{q,\CA}$.
Then, $[c\varphi ,x] = (q^m -1)xc\varphi + cd$.  Induction and the fact
that $(q-1)^n=0$ completes the lemma.
\end{proof}
\begin{remark}
In general, we can prove the following (\cite{I} lemma 1.0.0.23):  
Let $\K = \mathbb{Q}[t]/t^n$.  
Let $R$ be a $\Gamma$-graded \K -algebra.
Suppose $\beta :\Gamma \times \Gamma \to \K ^{\times}$
be a bicharacter such that $\beta (a ,b) - 1 \in (t)$.
Futher assume that $\bar{R} := R/t$ is 
commutative.
Then
\begin{align*}
D_q(R) = \{ \varphi \in \mathit{grHom}_{\K}(R,R)|
                 \textit{there exists an }
                n \textit{ such that }\\
                [\cdots [[\varphi ,a_1],a_2],\cdots ,a_n] =0 , a_i \in R \}
\end{align*}
\end{remark}
\begin{remark}
We have
\[
 \lim _{\leftarrow n}D_{q,\CA} ((A /(q-1)^n) [x]) \cong 
                \mathbb{Q}<x,\partial >[[t]] /[\partial  , x] =1,
\]
as $\mathbb{Q}$ -algebras.
\end{remark}
\begin{definition}\hfill
\begin{enumerate}
\item   For each $n\geq1$, let
        \begin{align*}
        \Gamma _{q,n} := \mathit{Ker} \hookrightarrow
                D_{q,\CA}(\CA /(q-1)^n [x]) &\bigoplus
                        D_{q,\CA}(\CA /(q-1)^n [x^{-1}]\\
                        &\to D_{q,\CA}(\CA /(q-1)^n [x,x^{-1}]);\\
                (\varphi _1,\varphi _2) &\mapsto \varphi _1 - \varphi _2,
        \end{align*}
        the algebra of $\CA$-linear global $q$-differential operators on
        $\mathbb{P}^1_{n\CA}$.
\item   The $\Gamma _{q,n}$ form an inverse system with 
        $\hat{\Gamma}$ as its inverse
        limit.
\end{enumerate}
\end{definition}
\begin{remark}\label{R:nton-1}\hfill
\begin{enumerate}
\item   If $(\varphi _1,\varphi _2) \in \Gamma _{q,n}$ then
        $\varphi _1 \in (\CA /(q-1)^n) <\partial _x,x^2\partial_x> $ and
        $\varphi _2 \in (\CA /(q-1)^n) 
        <\partial _{x^{-1}},x^{-2}\partial_{x^{-1}}> $.
\item   The map $\varphi \in (\CA /(q-1)^n)<\partial _x,x^2\partial _x>$
        if and only if $\varphi$ can be written as
        \[
        \varphi = \sum _{i\geq 1} f_i(q,q^{-1})g_i(x)\partial _x^i
        \]
        where $g_i(z)\partial _x^i \in \mathbb{Q}<\partial _x,x^2\partial_x>$.
        Moreover, if image of $\varphi$ is contained in 
        $(q-1)^{n-1}\CA /((q-1)^{n-1} [x]$, 
        then by induction on $n$, we can see that 
        $f_i = (q-1)^{n-1} h_i$.
        Thus, there exists $d\in \CA /(q-1)^n <\partial _x,x^2\partial _x>$
        such that $\varphi = (q-1)^n d$.
\end{enumerate}
\end{remark}
\subsection{Inverse limit of $U_q$}
Let $[m] = \frac{q^m-q^{-m}}{q-q^{-1}}$ and  
$[m]! = \prod _{1\leq i\leq m}[m]$. Denote by
\[
E^{(m)} = \frac{E^m}{[m]!}, F^{(m)} = \frac{F^m}{[m]!}, m \in \mathbb{Z}.
\]

Let $U_{q,\CA}$ be the subalgebra of $U_q$ generated by 
$E^{(m)},F^{(m)},K,K^{-1},m\in \mathbb{Z}$ over \CA.
For each $n\geq 1$,
let $U_{q,n}$ denote the ring $U_{q,\CA}/(q-1)^n$, and
$\hat{U_q}$ denote their inverse limit.
 There are homomorphisms
$\alpha _n, \gamma _n$ induced by $\alpha ,\gamma$ (defined in the subsection
\ref{S:UqtoDq}) respectively, giving a map
$\eta _n : U_{q,n} \to \Gamma _{q,n}$, whose inverse limit is
denoted by $\hat{\eta }: \hat{U_q} \to \hat{\Gamma} $.
\begin{theorem}\label{T:inverselimits}
The map
$\hat{\eta}: \hat{U_q} \to \hat{\Gamma}$ is a surjection.
\end{theorem}
\begin{proof}
We show by induction that $\eta _n$ is surjective for $n\geq 1$.
When $n=1$, note that $\Gamma _{q,1}$ is generated over $\mathbb{Q}$ by
$(\partial ,-x^{-2}\partial _{x^{-1}}), (-x^2\partial , \partial _{x^{-1}})$.
The map $\eta _1$ is clearly surjective.
Let $(\varphi _1,\varphi _2) \in \Gamma _{q,n}$.  Consider
$(\bar{\varphi _1},\bar{\varphi _2}) \in \Gamma _{q,n-1}$.
Since $\eta _{n-1}$ is surjective, there exists a $u\in U_{q,\CA}$
such that
$\eta _{n-1} (\bar{u}) = (\bar{\varphi _1},\bar{\varphi _2})$.
Consider $\eta _n (\bar{u}) - (\varphi _1,\varphi_2) =(\psi _1,\psi_2)
\in \Gamma _{q,n}$.
Since $(\bar{\psi_1},\bar{\psi _2}) =0 \in \Gamma _{q,n-1}$, we have
$(\psi_1,\psi_2) = (q-1)^{n-1}(d_1,d_2)$ for
$(d_1,d_2) \in \Gamma _{q,n}$ by remark \ref{R:nton-1}.
Consider $(\bar{d_1},\bar{d_2}) \in \Gamma _{q,n-1}$.
Again by surjectivity of $\eta _{n-1}$, there exists
a $v\in U_{q,\CA}$ such that
$\eta _{n-1} (\bar{v}) = (\bar{d_1},\bar{d_2})$.
Since $n\gneq 1$, we have $\eta _n ((q-1)^n \bar{v}) = (d_1,d_2)$.
Thus,
$\eta _n (\bar{u} - (q-1)^n\bar{v}) = (\varphi _1,\varphi _2)$.
\end{proof}

\bibliographystyle{amsplain}

\end{document}